\documentclass[letterpaper, 10 pt, conference]{ieeeconf}  

\IEEEoverridecommandlockouts                              
\usepackage{graphics} 
\usepackage{amsmath} 
\usepackage{amssymb}  
\usepackage{comment}
\usepackage{theorem}
\usepackage{tikz}
\usetikzlibrary{shapes,arrows,positioning}
\usepackage{tkz-graph}
\usepackage{booktabs}
\usepackage{hyperref}  
\usepackage{url}

\newcommand{\bea}{\begin{eqnarray}}
\newcommand{\eea}{\end{eqnarray}}
\newcommand{\beas}{\begin{eqnarray*}}
\newcommand{\eeas}{\end{eqnarray*}}
\newcommand{\leftm}{\left[\begin{array}}
\newcommand{\rightm}{\end{array}\right]}

\newcommand{\mL}{\mathcal{L}}

\newcommand{\be}[1]{\begin{equation}\label{#1}}
\newcommand{\benon}{\begin{equation*}}  
\newcommand{\bemuln}[1]{\begin{multline}\label{#1}}
\newcommand{\bemul}{\begin{multline*}}
\newcommand{\bee}{\begin{eqnarray*}}
\newcommand{\eee}{\end{eqnarray*}}
\newcommand{\been}[1]{\begin{eqnarray}\label{#1}}
\newcommand{\eeen}{\end{eqnarray}}
\newcommand{\began}[1]{\begin{gather}\label{#1}}
\newcommand{\bega}{\begin{gather*}}
\newcommand{\bealn}[1]{\begin{align}\label{#1}}
\newcommand{\beal}{\begin{align*}}
\newcommand{\bealatn}[2]{\begin{alignat}{#1}\label{#2}}
\newcommand{\bealat}{\begin{alignat*}}
\newcommand{\bexalatn}[1]{\begin{xalignat}\label{#1}}
\newcommand{\bexalat}{\begin{xalignat*}}

\newcommand{\Ra}{\Rightarrow}

\newcommand{\mb}{\mathbf}

\newcommand{\mbb}{\mathbb}

{\theoremstyle{plain}

}
{\theoremstyle{break} \theorembodyfont{\it}

 }

\def\bg{{\mathbf g}}
\def\bh{{\mathbf h}}

\def\bt{{\mathbf t}}

\def\bw{{\mathbf w}}
\def\bx{{\mathbf x}}  
\def\by{{\mathbf y}}
\def\bz{{\mathbf z}}
\def\bA{{\mathbf A}}
\def\bB{{\mathbf B}}

\def\bI{{\mathbf I}}

\def\bN{{\mathbf N}}

\def\texitem#1{\par\smallskip\noindent\hangindent 25pt
               \hbox to 25pt {\hss #1 ~}\ignorespaces}

\newcommand{\bzero}{{\mathbf{0}}}

\newcommand{\scrA}{\mathcal{A}}

\newcommand{\scrD}{\mathcal{D}}

\newcommand{\scrF}{\mathcal{F}}

\newcommand{\scrH}{\mathcal{H}}

\newcommand{\scrR}{\mathcal{R}}

\newcommand{\scrW}{\mathcal{W}}

\newcommand{\bbeta}{\boldsymbol{\beta}}

\newcommand{\bepsilon}{\boldsymbol{\epsilon}}

\newcommand{\bnu}{{\boldsymbol{\nu}}}

\title{\LARGE \bf Joint Estimation of OD Demands and Cost Functions in Transportation
  Networks from Data~\authorrefmark{1} 
\thanks{* Research partially supported by the
    NSF under grants DMS-1664644 and CNS-1645681, by the ONR under MURI grant
    N00014-16-1-2832, and by the Boston University Division of Systems Engineering.}}

\author{Salom\'{o}n Wollenstein-Betech$^{1}$, Chuangchuang Sun$^2$, Jing Zhang$^3$, and Ioannis Ch. Paschalidis$^4$ 
\thanks{$^{1,2}$ Division of Systems Engineering, Boston University, {\tt
    \{salomonw, ccsun\}@bu.edu}.} 
\thanks{$^{3}$ Mitsubishi Electric Research Laboratories, Cambridge, MA,  
{\tt jingzhang@merl.com}.}
\thanks{$^{4}$Dept. of Electrical and Computer Engineering,
	Division of Systems Engineering,
	and Dept. of Biomedical Engineering, Boston University, 8 St. Mary's St.,
	Boston, MA 02215, USA.
	{\tt yannisp@bu.edu}, \url{http://sites.bu.edu/paschalidis/}.}}

\begin{document}

\maketitle
\thispagestyle{empty}
\pagestyle{empty}

\begin{abstract}
Existing work has tackled the problem of estimating Origin-Destination (OD) demands
and recovering travel latency functions in transportation networks under the
\emph{Wardropian} assumption. The ultimate objective is to derive an accurate
predictive model of the network to enable optimization and control. However, these
two problems are typically treated separately and estimation is based on parametric
models. In this paper, we propose a method to jointly recover nonparametric travel
latency cost functions and estimate OD demands using traffic flow data. We formulate
the problem as a bilevel optimization problem and develop an iterative first-order
optimization algorithm to solve it. A numerical example using the Braess Network is
presented to demonstrate the effectiveness of our method.
\end{abstract}


\section{INTRODUCTION}  \label{sec:intro}
The purpose of solving the \emph{Traffic Assignment Problem} (TAP) in transportation
planning processes is to evaluate performance metrics of the system, assess
deficiencies and evaluate potential improvements and capacity expansions to the
transportation network.

The TAP assumes that users selfishly choose the best route in the network resulting
in an equilibrium known as \emph{Wardrop equilibrium}. Modeling drivers' routing
behavior under the Wardrop equilibrium assumption is one of the most widely-used
frameworks for the purpose of analyzing transportation networks, with applications in
traffic diagnosis, control, and optimization \cite{Merchant78,Patriksson1994}.  This
modeling framework uses three main inputs: $(1)$ a strongly connected directed graph;
$(2)$ an Origin Destination (OD) traffic demand vector; and $(3)$ a link
\emph{latency cost} or \emph{travel time cost} function that typically depends on
link flows. Small perturbations to these OD demand estimates and \emph{travel time
  functions} may have a large impact on the equilibrium solution
\cite{bertsimas2015data}.

In practice, however, OD demands and cost functions are not readily available. The OD
demand estimation problem for the static TAP has been solved differently depending on
whether a network is congested or not. For uncongested networks, entropy maximization
\cite{VanZuylen80}, generalized least squares \cite{Hazelton2000} and maximum
likelihood estimation \cite{Spiess87} have been used. Whereas for congested networks,
estimating OD demands has been done by solving a bilevel optimization problem given
the circular dependence between the OD estimation and the traffic flow assignment
\cite{Daamen2014}.

The problem of estimating \emph{travel time functions} has received less attention in
the transportation community. In the context of transportation systems, as traffic
volume grows we expect the speed on the link to decrease, first slowly but as queues
start to accumulate, the effects become more significant. Therefore, these functions
are usually modeled as positive, nonlinear and strictly increasing functions.
A typical \emph{travel time function} is as a polynomial function. In particular, urban planners and researchers often use the {\em Bureau of Public
Roads (BPR)} function \cite{BPR}:
\begin{equation} \label{bpr}
t(x_a) = t_a^0(1+0.15(x_a/m_a    )^4),
\end{equation}
where $t_a^0$ is the free-flow travel time, $x_a$ the flow, and $m_a$ the capacity of
link $a$.

With the increasing availability of various sensors, large traffic datasets have been
collected, raising the possibility of estimating OD demands and \emph{travel time
  functions} from data by solving appropriate inverse optimization problems.  More
specifically, given an OD demand and equilibrium flows, recovering the \emph{travel
  time function} can be performed for both single-class vehicle networks
\cite{bertsimas2015data,Zhang2016} and multi-class vehicle networks \cite{Zhang2018}.

Most of the existing work typically deals with these two inverse problems separately;
a limitation we seek to address in this paper.  Closer to the goal of our work,
\cite{Yang2001} considered the simultaneous estimation of travel cost and OD demand
in a \emph{Stochastic User Equilibrium} setting. Yet, this work does not attempt to
estimate (nonparametrically) the full structure of the travel cost functions as we
do. Rather, it seeks to estimate a sensitivity constant that adjusts how a given
travel cost function affects route choice probabilities.

In this paper, we aim to jointly investigate the two related inverse problems --
recovering cost functions (IP-1) in a non-parametric setting and adjusting OD demand
matrices (IP-2). Our work contributes to improving the consistency and robustness of
the data-driven traffic model. The ultimate utility of obtaining such a model is to
use it to make predictions under various topology and demand scenarios, drive control
and optimization tasks, or simply assess the amount of inefficiency of the system
(e.g., as in ~\cite{Zhang2018}). In this work we consider only the (data-driven)
model estimation problem.

We solve the joint problem by converting the bilevel optimization model into a
single-level one. We do this by transforming the lower-level problem (IP-1) into
constraints for the upper-level one (IP-2). As a result, we obtain a formulation with
a quadratic objective and non-convex constraints. Using weak duality and an iterative
approach, we are able to relax the non-convex constraints, which allows the problem
to be solved using a first-order feasible direction algorithm.  To validate its
effectiveness and performance, we conduct a numerical experiment using the Braess'
network \cite{Braess2005}. In this example, we show that the algorithm approaches the
\emph{ground truth} values of both \emph{travel time functions} and OD demands.

The rest of the paper is organized as follows. In Sec.~\ref{sec:model and pre} we
introduce the modeling framework and mathematical definitions used throughout the
paper. In Sec.~\ref{sec:JointProblem} we present the structure of the joint problem,
its transformation to its Frank-Wolfe form, and a method for calculating the gradient
of the cost function. In Sec.~\ref{sec:numericalExample} we present some numerical
results applied to the Braess network. Conclusions are in Sec.~\ref{sec:conclusion}.

\textbf{Notation:} All vectors are column vectors and denoted by bold lowercase
letters. Bold uppercase letters denote matrices. To economize space, we write $\bx =
(x_1, \ldots ,x_{\text{dim}(\bx)})$ to denote the column vector $\mathbf{x}$, where
$\text{dim}(\bx)$ is its dimensionality.  We use ``prime'' to denote the transpose of
a matrix or vector. We denote by $\bzero$ and $\bI$ the vector of all zeroes and the
identity matrix, respectively. Unless otherwise specified, $\|\cdot\|$ denotes the $\ell_2$
norm. $|\scrD|$ denotes the cardinality of a set $\scrD$, and
$\left[\kern-0.15em\left[ \mathcal{D} \right]\kern-0.15em\right]$ the set $\{ 1,
\ldots , |\mathcal{D}| \}$.

\section{MODEL AND PRELIMINARIES} \label{sec:model and pre}
\subsection{Transportation network model and definitions} \label{subsec:network model}

Consider a strongly-connected directed graph denoted by $G\left( \mathcal{V},
\mathcal{A} \right)$, where $\mathcal{V}$ is the set of nodes and $\mathcal{A}$ is
the set of links.  Let $\mathbf{N} \in {\left\{ {0,1, - 1} \right\}^{\left|
    \mathcal{V} \right| \times \left| \mathcal{A}\right|}}$ be the node-link
incidence matrix, and let $\textbf{e}_{a}\in \mbb{R}^{|\scrA|}$ be a vector with an
entry equal to $1$ corresponding to link $a$ and all the other entries set to
$0$. Let $\mathbf{w} = (w_s,w_t)$ denote an Origin-Destination (OD) pair and
$\mathcal{W} = \left\{ {{\mathbf{w}_i}:{\mathbf{w}_i} = \left( {{w_{si}},{w_{ti}}}
  \right), \,i \in \left[\kern-0.15em\left[ \mathcal{W} \right]\kern-0.15em\right]}
\right\}$ be the set of all OD pairs.  Furthermore, let ${d^{\mathbf{w}}} \ge 0$ be
the flow demand that travels from origin $w_s$ to destination $w_t$. In the same
manner, let us denote by $\mathbf{d^w} \in \mathbb{R}^{|\mathcal{V}|}$ the vector of
all zeros except for the coordinates of nodes $w_s$ and $w_t$ which take values
$-d^{\mathbf{w}}$ and $d^{\mathbf{w}}$, respectively.  We will also use vector
$\bg=(d^{\bw}; \bw\in \scrW)$ to denote the flow demands for all OD pairs. Let $x_a$
be the total link flow of link $a \in \mathcal{A}$ and $\mathbf{x}$ the vector of
these flows.  Let $\mathcal{F}$ be the set of feasible flow vectors defined by \\
\[
\mathcal{F} = \Big\{  \mathbf{x} \in \mathbb{R}_+^{|\mathcal{A}|}: \mathbf{x} =
		\sum\limits_{\mathbf{w} \in \mathcal{W}} {{\mathbf{x}^{\mathbf{w}}}},\, 
		\mathbf{N}{\mathbf{x}^{\mathbf{w}}} = {\mathbf{d}^{\mathbf{w}}},\, \forall \mathbf{w} \in \mathcal{W} \Big\},
\]
where $\mathbf{x}^\mathbf{w}$ is the flow vector attributed to OD pair $\mathbf{w}$.

For each OD pair $\bw$ let us also define a set of possible routes
$\mathcal{R}^{\bw}$; each route $r \in \mathcal{R}^{\bw}$ is a sequence of links
starting from the origin $w_s$ and ending at the destination $w_t$. We will write
$a\in r$ if a route $r$ contains link $a$. For each OD pair $\bw_i\in \scrW$ we
define the indicator functions
\begin{equation} \label{delta}
   \delta_{r}^{ai} = \begin{cases}
        1, & \text{if }  r \in \mathcal{R}^{\bw_i} \text{ uses link } a\\
        0, & \text{otherwise.} 
\end{cases}
\end{equation}

Finally, we denote with $t_a(\mathbf{x}) : \mathbb{R}_{+}^{\left| \mathcal{A}
  \right|} \mapsto \mathbb{R}_{+}$ the \emph{latency cost} (i.e., travel time)
function for link $a$ and write $\bt(\cdot)$ for the vector of these link functions.
Using the same structure used in \cite{Beckmann1955} we can
characterize $t_a(x_a)$ as:
\[
\label{eq:latency function}
t_a(x_a) = t_a^0f(x_a/m_a),
\]
where $m_a$ is the flow capacity of link $a$, $f(\cdot)$ is a strictly increasing,
positive, and continuously differentiable function, and $t_a^0$ is the free-flow
travel time on link $a$. We set $f(0)=1$, which ensures that if there is no
constraint on flow capacity, the travel time $t_a$ is equal to the free-flow travel
time.

\subsection{Wardrop equilibrium} \label{subsec: wardprop}

The notion of a Wardrop equilibrium, sometimes referred to as a non-atomic
game\footnote{These are games where every user (driver) has a negligible
  contribution to the overall traffic. Hence, the actions of individual users have
  essentially no effect on network congestion.}, is interpreted as requiring that all
users optimize their travel times. In general, a feasible flow
$\mathbf{x^{*}}$ is a Wardrop equilibrium if for every OD pair $\mathbf{w}_i$, and
any route $r\in \scrR^{\bw_i}$ with positive flow, the latency cost (i.e., travel
time) is no greater than the travel time on any other route. It is worth mentioning
that given $G(\mathcal{V},\mathcal{A})$ and $f(\cdot)$ there exists a unique 
equilibrium\footnote{Backman proves this using KKT conditions
  \cite{Beckmann1955}.}. Such a result is the solution to the \emph{Traffic Assignment
  Problem (TAP)} which precisely returns the flows that minimize the \emph{potential
  function}:
\[
\Phi(\mathbf{x}) = \sum\limits_{a \in \mathcal{A}} \int\limits_{0}^{x_a}t_a(s)ds,
\]
where the integral is adding the costs of the flow segments of link $a$.  The
function $f(\cdot)$ is continuous and $\mathcal{F}$ is a compact set, thus,
Weierstrass Theorem implies there exists a solution. Moreover, since cost functions
are non-decreasing (by assumption), then $\Phi(\cdot)$ is convex and therefore a
unique solution exists~\cite{Beckmann1955}.

\subsection{Models} \label{subsec:models}

\subsubsection{User-centric} \label{subsubsec:user-centric}
As stated in the previous section, the TAP (also known as the \emph{user-centric
  forward optimization problem}) can be formulated as 
\begin{equation}\label{TAP}
\min_{\mathbf{x} \in \mathcal{F}} \ \ \sum\limits_{a \in \mathcal{A}}{\int\limits_{0}^{x_a}{t_a(s)ds}}.
\end{equation}
An alternative way of solving this problem is via a \emph{Variational
  Inequality} (VI) formulation as first proposed in \cite{Smith1979,Dafermos1980};
finding a solution $\bx^*$ to  
\begin{equation} \label{VI}
\mathbf{t}(\mathbf{x}^{*})'(\mathbf{x}-\mathbf{x^*}) \geq 0, \ \ \forall \mathbf{x}
\in \mathcal{F}. 
\end{equation} 
In order for the solution of (\ref{VI}) to be equivalent to the solution of
(\ref{TAP}) we have to assume $(i)$ strong monotonicity of $\mathbf{t}(\cdot)$ over
$\mathcal{F}$, $(ii)$ $\mathbf{t}(\cdot)$ to be continuously differentiable over
$\mathbb{R}_{+}^{|\mathcal{A}|}$, and $(iii)$ $\scrF$ to contain an interior point
(Slater's condition). One of the most successful algorithms to find such an equilibrium is the
\textit{Method of Successive Averages (MSA)} proposed in \cite{Daganzo77} which uses
a Frank-Wolfe type algorithm.

\subsubsection{User-Centric Inverse Model (I-VI)} \label{subsubsec:  I-VI}

Given that one of the parameters of the TAP is the latency cost functions, we
aim to estimate them (in particular function $f(\cdot)$) using data. To that end, we
consider an \emph{Inverse Variational Inequality} problem (I-VI).  We assume that the
data measurements are solutions of the TAP for specific cost functions and OD
demands.  Therefore, it is natural to think about these flows as snapshots of the
network at different instants.  Let $k \in \left[\kern-0.15em\left[\mathcal{K}
    \right]\kern-0.15em\right]$ index different snapshots of a network with
corresponding flows $\mathbf{x}^{(k)}=(x_a^{(k)};\ a \in \scrA^{(k)})$, where the set
$\scrA^{(k)}\subset \scrA$ denotes the links on which we have flow measurements for
instance $k$. (We will use $\scrF^{(k)}$, $\bN_k$, and $\scrW^{(k)}$ to denote the
set of feasible flows, node-link incidence matrix, and OD pairs for the network
instance $k$.)  The inverse formulation of the \textit{Wardrop equilibrium} seeks to
find a cost function $\bt(\cdot)$ (or, equivalently, $f(\cdot)$) such that each flow
observation is as close to an equilibrium as possible.  Because this formulation
relies on measured data, we expect measurement noise. Hence, the notion of an
approximate solution to this problem is natural. For a given $\epsilon>0$, we define
an $\epsilon$-approximate solution $\hat{\bx}$ to the VI as satisfying:
\begin{equation} \label{VI-eapprox}
\mathbf{t}(\hat{\mathbf{x}})'(\mathbf{x}-\hat{\mathbf{x}}) \geq -\epsilon, \ \ \forall \mathbf{x} \in \mathcal{F}.
\end{equation}

The inverse VI problem amounts to finding a function ${f(\cdot)}$ such that
$\mathbf{x}^{(k)}$ is an $\epsilon_k$-approximate solution to VI$(\mathbf{t},
\mathcal{F}^{{(k)}})$ for each $k$.  Denoting $\boldsymbol{\epsilon}
\stackrel{\triangle}{=} (\epsilon_k; \, k \in [\kern-0.15em[ \mathcal{K}
  ]\kern-0.15em])$, we can formulate the inverse VI problem as in
\cite{bertsimas2015data,Zhang2018}. Then we define the (I-VI) problem as minimizing the
$\ell_2$ norm of $\boldsymbol{\epsilon}$:
\begin{align} \label{I-VI-eapprox}
\min_{\bt(\cdot),\boldsymbol{\epsilon}}  & \ \  \|\boldsymbol{\epsilon} \| \\ \notag
\text{s.t. } & \ \  \mathbf{t}(\mathbf{x}^{(k)})'(\mathbf{x}-\mathbf{x}^{(k)}) \geq -\epsilon_k,  \hspace{5mm} \forall \mathbf{x} \in \mathcal{F}^{(k)}, k \in \left[\kern-0.15em\left[\mathcal{K} \right]\kern-0.15em\right], \\
& \ \ \epsilon_k > 0, \qquad \forall k \in \left[\kern-0.15em\left[\mathcal{K}
    \right]\kern-0.15em\right]. \notag 
\end{align}
Notice that in this formulation, the set of constraints restricts the travel time
function to be within $\epsilon_k$ units of the \emph{Wardrop equilibrium} flows for
each sample. In this sense, if we solve the problem using $k$ of these constraints
for multiple observed networks, we will find a more ``stable'' travel time function.

In order to solve this problem we express the function $f(\cdot)$ in a
\emph{Reproducing Kernel Hilbert Space} (RKHS) $\mathcal{H}$ as in
\cite{bertsimas2015data}. This leads to the following formulation of the
$\epsilon$-approximate Inverse Variational Inequality Problem ($\epsilon$I-VI):
\begin{align} \label{invIV-I} 
\min_{f,\by,\bepsilon} & \ \  \| \bepsilon \| + \gamma \| f \|^{2}_{\scrH} \\
\text{s.t. }  & \ \ \mathbf{e}_a' \mathbf{N}_k' \mathbf{y}^{\mathbf{w}} \leq t^0_a
f\left(\frac{x_a}{m_a}\right), \forall \mathbf{w} \in \mathcal{W}^{(k)}, a \in
\mathcal{A}^{(k)}, k, \notag \\
& \ \ \sum_{a \in \mathcal{A}^{(k)}}  t^0_a x_a f\left(\frac{x_a}{m_a}\right)
- \hspace{-0.3cm} \sum_{\bw \in
  \mathcal{W}^{(k)}} (\mathbf{d^w})'\mathbf{y^w} \leq \epsilon_k, \forall k, \notag \\
& \ \ f\left(\frac{x_a}{m_a}\right) \leq  f\left(\frac{x_{\hat{a}}}{m_{\hat{a}}}\right), \forall a, \hat{a}
\in \cup_{k} \mathcal{A}^{(k)}\, \text{s.t.}\, \frac{x_a}
    {m_a}\leq\frac{x_{\hat{a}}}{m_{\hat{a}}},  
\notag \\
& \ \ \boldsymbol{\epsilon} \geq 0,\;  f \in \mathcal{H}, f(0) = 1, \notag
\end{align}
where the first constraint corresponds to dual feasibility, the second constraint
maintains the primal-dual gap within $\epsilon$, and the third constraint imposes the
assumption that $f(\cdot)$ is monotone. We note that $\mathbf{y^w}$ contains dual
variables associated with the VI problem, $\|\cdot\|_{\scrH}$ is the norm of the
RKHS, and $\gamma$ is a regularization parameter. A larger $\gamma$ will recover a
more general $f(\cdot)$ whereas a smaller one will recover an $f(\cdot)$ which fits
the dataset better.

As we can see, the problem we have defined is still hard to solve since it involves
optimization over functions $f(\cdot)$. However, we specify $\mathcal{H}$ (and thus
the class of $f(\cdot)$) by choosing a polynomial kernel \cite{bertsimas2015data}, i.e.,
using kernel functions $\phi(x,y) = (c+xy)^n$. We believe this is a good choice 
since it matches our intuition on how congestion affects the latency cost of links
(cf. (\ref{bpr})). The polynomial kernel function can be rewritten as
\[
\phi(x,y) = (c+xy)^n = \sum\limits_{i=0}^{n} {n\choose i}c^{n-1}x^iy^i.
\]
Then, using the representer theorem for kernel functions, we can modify the cost
function of the ($\boldsymbol{\epsilon}$I-VI) problem to a quadratic function
parameterized by $\boldsymbol{\beta}= \{\beta_j : j=1,\ldots,n \}$
resulting in a tractable Quadratic Programming (QP) problem (see
\cite{bertsimas2015data,Zhang2018} for details). As an output to this reformulated
($\boldsymbol{\epsilon}$I-VI) problem we obtain
$\boldsymbol{\beta}^{*}$, and therefore our estimator for $f(\cdot)$ is
equal to
\[
\hat{f}(x) = \sum\limits_{i=0}^{n}\beta_i^{*}x^{i} = 1 + \sum\limits_{i=1}^{n}\beta_{i}^{*}x^{i},
\]
where we set $\beta_0 = 1$ to have $f(0)=1$.

To facilitate the analysis of the joint problem presented in the next section, let us
write the QP problem corresponding to ($\boldsymbol{\epsilon}$I-VI) using compact
notation:
\begin{align} \label{compact} 
\min_{\boldsymbol{\beta},\by,\boldsymbol{\epsilon}} & \ \  \boldsymbol{\epsilon}'\mathbf{I}\boldsymbol{\epsilon} + \boldsymbol{\beta}'\mathbf{H}\boldsymbol{\beta} \\
\text{s.t. } & \ \  \mathbf{A}(\bg) \by + \mathbf{B} (\bx) \boldsymbol{\beta}+ \mathbf{C}\boldsymbol{\epsilon} + \bh \leq \boldsymbol{0}, \notag
\end{align}
where matrices $\mathbf{A}(\bg)$ and $\mathbf{B}(\bx)$ depend on the OD demand vector
$\bg$ and the provided data flow measurements $\mathbf{x}$, respectively, and
$\mathbf{H}$ is a positive definite matrix. We call this problem  (IP-1).

\section{THE JOINT PROBLEM} 
\label{sec:JointProblem}

\subsection{Bilevel formulation}
Unlike previous work, we will jointly recover both the \emph{travel time} function
$f(\cdot)$, specifically the coefficients $\boldsymbol{\beta} = (\beta_o, \ldots,
\beta_n)$, and the OD demand vector $\bg$. To simplify notation, we let
$\mathbf{x}(\boldsymbol{\beta}, \bg) = (x_a(\boldsymbol{\beta},
\bg);\ \forall a \in \mathcal{A})$ be the optimal solution to the
VI$(\mathbf{t}, \mathcal{F})$ (i.e., the TAP), for any given feasible
$\boldsymbol{\beta}$ and $\bg$. Recall that we observe an equilibrium flow
vector from data which we define as $\mathbf{x^*}=(x_a^*;\ \forall a \in
\mathcal{A})$. Equipped with these definitions we can define the bilevel optimization
problem as follows
\begin{align}\label{eq:bilevel}
\min_{\bbeta,\bg}& \ \ F(\boldsymbol{\beta},
\bg) \stackrel{\triangle}{=}    \sum_{a\in \scrA}
(x_{a}(\boldsymbol{\boldsymbol{\beta}},\bg) - x_{a}^*)^2  \\ 
\text{s.t. } &  (\boldsymbol{\beta},\by,\boldsymbol{\epsilon}) = \arg\min_{\boldsymbol{\beta},\by,\boldsymbol{\epsilon}}\big\{\boldsymbol{\epsilon}'\mathbf{I}\boldsymbol{\epsilon} + \boldsymbol{\beta}'\mathbf{H}\boldsymbol{\beta}, \notag \\ 
& \text{s.t. } \mathbf{A}(\bg)\by + \mathbf{B}(\bx(\bbeta, \bg)) \boldsymbol{\beta}+ \mathbf{C}\boldsymbol{\epsilon} + \bh \leq \boldsymbol{0} \big\}, \notag\\
&  \boldsymbol{\boldsymbol{\beta}} \ge \bold{0},\  \bg \ge \bold{0}. \notag
\end{align}
Notice that $F(\boldsymbol{\beta}, \bg)$ is bounded below by $0$. 

To solve this problem we replace the convex lower-level problem (IP-1) by its KKT
optimality conditions and write the bilevel problem as a single-level
problem. Finally, we relax the resulting formulation to make it solvable by using a
feasible direction method (Frank-Wolfe).

\subsection{IP-1 Optimality conditions}

To reduce the lower level problem in \eqref{eq:bilevel} into its equivalent
optimality conditions, we first write the Lagrangian function:
\begin{equation*} \label{eq:Lagrangian}
\mL(\boldsymbol{\beta},\by,\boldsymbol{\epsilon}; \boldsymbol{\nu}) =  \boldsymbol{\epsilon}'\mathbf{I}\boldsymbol{\epsilon} + \boldsymbol{\beta}'\mathbf{H}\boldsymbol{\beta}  
+  \boldsymbol{\nu}'(\mathbf{A}\by + \mathbf{B}\boldsymbol{\beta}+ \mathbf{C}\boldsymbol{\epsilon} + \bh), 
\end{equation*}
where $\boldsymbol{\nu}$ are the dual variables and (for ease of notation) we dropped
the dependence of $\bA$ and $\bB$ on $\bg$ and $\bx(\bbeta, \bg)$, respectively. 

This leads to the first order optimality conditions:
\bea\label{eq:KKT}
\partial\mL/\partial\boldsymbol{\epsilon} &=& 2\mathbf{I} \boldsymbol{\epsilon} + \mathbf{C}'\boldsymbol{\nu} = \mathbf{0} \Ra  \boldsymbol{\epsilon} = -(1/2)\mathbf{I}^{-1}\mathbf{C}' \boldsymbol{\nu}, \nonumber\\
\partial\mL/\partial\boldsymbol{\beta} &=& 2\mathbf{H}\boldsymbol{\beta} + \mathbf{B}' \boldsymbol{\nu} = \mathbf{0} \Ra  \boldsymbol{\beta} = -(1/2)\mathbf{H}^{-1}\mathbf{B}'\boldsymbol{\nu}, \nonumber\\
\partial\mL/\partial \by &=& \mathbf{A}'\boldsymbol{\nu} = \mathbf{0}.
\eea
Substituting $\boldsymbol{\beta}$ and $\boldsymbol{\epsilon}$ in the Lagrangian using
\eqref{eq:KKT}, we can write the dual objective function as 
\bea\label{eq:dual_objective} 
D(\boldsymbol{\nu}) = -\frac{1}{4}\boldsymbol{\nu}'\mathbf{CIC}' \boldsymbol{\nu} - \frac{1}{4}\boldsymbol{\nu}'\mathbf{BH^{-1}B} \boldsymbol{\nu} + \bh'\bnu. 
\eea

Consequently, for each primal-dual pair $(\boldsymbol{\beta},\by,\boldsymbol{\epsilon};\, \boldsymbol{\nu})$ in the lower-level optimization problem, it is sufficient and necessary to satisfy the  conditions
\bea\label{eq:kkt_cond}
&\mathbf{A}\by + \mathbf{B}\boldsymbol{\beta}+ \mathbf{C}\boldsymbol{\epsilon} + \bh \leq \boldsymbol{0},  \\
&\mathbf{A}'\boldsymbol{\nu} = \boldsymbol{0}, \notag \\
&\boldsymbol{\nu} \ge \mathbf{0}, \notag\\
&\boldsymbol{\epsilon}'\mathbf{I}\boldsymbol{\epsilon} +
\boldsymbol{\beta}'\mathbf{H}\boldsymbol{\beta} \hspace{-0.4pt} = \hspace{-0.4pt} -\frac{1}{4}\boldsymbol{\nu}'\mathbf{CIC}'\boldsymbol{\nu} -\frac{1}{4}\boldsymbol{\nu}'\mathbf{BH^{-1}B}'\boldsymbol{\nu} + \bh' \boldsymbol{\nu}, \notag
\eea
to reach optimality.

\subsection{Relaxation and Frank-Wolfe}

So far, we have eliminated the lower optimization problem by transforming it into
constraints involving the dual variables. Note that the fourth constraint of
\eqref{eq:kkt_cond}, corresponding to strong duality of (IP-1), is a non-convex
quadratic equality constraint. To address this issue, we relax it by requiring that
the duality gap is upper bounded by some $\xi$ and penalizing $\xi$:
\begin{align} \label{eq:bilevel_with_kkt}
\min_{\bbeta, \by, \bepsilon, \bg, \bnu, \xi}& F(\boldsymbol{\beta}, \bg, \xi)
\stackrel{\triangle}{=}    \sum_{a\in \scrA} (x_{a}(\bbeta,\bg) - x_{a}^*)^2 + \lambda \xi  \\
\text{s.t. } &\mathbf{A}\by + \mathbf{B} \boldsymbol{\beta}+ \mathbf{C}\boldsymbol{\epsilon} + \bh \leq \boldsymbol{0}, \notag\\
&\mathbf{A}'\boldsymbol{\nu} = \boldsymbol{0}, \notag\\
&\boldsymbol{\epsilon}'\mathbf{I}\boldsymbol{\epsilon} +
\boldsymbol{\beta}'\mathbf{H}\boldsymbol{\beta} +\frac{1}{4}
\boldsymbol{\nu}'\mathbf{CIC}'\boldsymbol{\nu} \notag\\
& + \frac{1}{4}
\boldsymbol{\nu}'\mathbf{B H^{-1}B}'\boldsymbol{\nu} - \bh'\bnu \leq \xi, \notag \\
&\boldsymbol{\nu}, \boldsymbol{g}, \boldsymbol{\beta}, \xi \geq
\mathbf{0}, \notag
\end{align}
where, again, we have suppressed the dependence of $\bA$ and $\bB$ on $\bg$ and
$\bx(\bbeta, \bg)$, respectively. Notice that both the objective and the constraints
(through $\bA$ and $\bB$) are nonlinear functions of $\bbeta, \bg$ through
$\bx(\bbeta, \bg)$.

We next develop an iterative \emph{feasible direction} method. Let
$\bz=(\bbeta,\bg,\xi)$ and $j$ denote the iteration count. 
We evaluate
the gradient of $F(\cdot)$ at the previous iteration and seek the steepest feasible
direction of descent by solving: 
\begin{align} \label{eq:FrankWolfe}
\min_{\bz_j, \by, \boldsymbol{\nu}, \boldsymbol{\epsilon}} & \ \nabla F(\bz_{j-1})' (\bz_{j-1}-\bz_{j}) \\
\text{s.t. }&  \mathbf{A} \by + \mathbf{B}\boldsymbol{\beta}+ \mathbf{C}\boldsymbol{\epsilon} + \bh \leq \boldsymbol{0}, \notag\\
&\mathbf{A}'\boldsymbol{\nu} = \boldsymbol{0}, \notag\\
&\boldsymbol{\epsilon}'\mathbf{I}\boldsymbol{\epsilon} +
\boldsymbol{\beta_j}'\mathbf{H}\boldsymbol{\beta_j}  +\frac{1}{4}
\boldsymbol{\nu}'\mathbf{CIC}'\boldsymbol{\nu}  \notag\\
& + \frac{1}{4} \boldsymbol{\nu}'\mathbf{BH^{-1}B}'\boldsymbol{\nu} 
- \bh'\bnu \leq \xi_{j} \notag \\
& \bg_{j-1} - c_1 \mb{e} \leq \bg_j \leq \bg_{j-1} +
c_2 \mb{e} \notag\\
&\boldsymbol{\nu}, \bz_j \geq \mathbf{0}, \notag
\end{align}
where we use $\mb{e}$ to denote the vector of all ones, $c_1, c_2$ are constants,
$\bA$ and $\bB$ in the constraints of (\ref{eq:FrankWolfe}) are functions of
$(\bbeta, \bg)$ evaluated at $(\bbeta_{j-1}, \bg_{j-1})$, and 
\begin{small}\begin{multline} \label{gradient}
\nabla F(\bz_j)' =
\bigg[\sum_{a \in \mathcal{A}} 2 (x_a(\bz_j)-x_a^*)
 \frac{\partial x_a (\bbeta_j,\bg_j)}{\partial \beta_{l}},\, l=1,\ldots,n;\\ 
\sum_{a \in \mathcal{A}}  2 (x_a(\bz_j)-x_a^*) \frac{\partial x_a
     (\bbeta_j,\bg_j)}{\partial g_{i}},\, i=1,\ldots,|\scrW|; \lambda \bigg].
\end{multline}
\end{small}
As a result, problem (\ref{eq:FrankWolfe}) has a linear objective and constraints
that are linear and convex quadratic, rendering it easy to solve. Given these
``constant'' approximations of the constraints at the prior iterate, the role of
$c_1, c_2$ is to ensure that the optimization takes place in a relatively small
``trust'' region for $\bg_j$ that is not too far from the prior iterate $\bg_{j-1}$.

\subsection{Derivatives}
For the cost function of \eqref{eq:FrankWolfe} (cf. (\ref{gradient})) we need to
estimate the partial derivatives of the link flows with respect to parameters
$\boldsymbol{\beta}$ of the latency functions and the OD demand vector $\bg$.

\subsubsection{Directional flow derivatives with respect to perturbations in OD
  demand}

Let us first derive an approximation to the gradient of
$\mathbf{x}(\boldsymbol{\beta},\bg)$ with respect to $\bg$. By adding the flows of
different OD pairs demands we have
\begin{align*}\label{linkflow}
x_a(\boldsymbol{\beta},\bg) & = \sum_{\{i: \bw_i \in \mathcal{W}\}} \sum_{r \in
  \mathcal{R}^{\bw_i}} \delta_{r}^{ai} p^{ir} g_{i} \notag \\
  & = \sum_{\{i: \bw_i \in \mathcal{W}\}} g_{i} \sum_{r \in
  \mathcal{R}^{\bw_i}} \delta_{r}^{ai} p^{ir},  
\end{align*}
where $\mathcal{R}^{\bw_i}$ denotes the set of feasible routes associated with OD
pair $\bw_i$, $\delta_{r}^{ai}$ was defined in (\ref{delta}), and $p^{ir}$ is the
probability that commuter in OD pair $\bw_i$ selects route $r \in
\mathcal{R}^{\bw_i}$.
	
For each OD pair $\bw_i \in \mathcal{W}$, let us only use the shortest route
$r_{i}(\boldsymbol{\beta},\bg)$ based on the travel latency cost (i.e., travel
time). Then we have  
\[
	\frac{{\partial {x_a}\left(\boldsymbol{\beta}, \bg  \right)}}{{\partial {g_{i}}}} \approx {\delta _{{r_{i}(\boldsymbol{\beta},\bg)}}^{ai}} = \begin{cases}
	1, & \text{if } a \in r_{i}(\boldsymbol{\beta},\bg), \\
	0, & \text{otherwise,} 
	\end{cases}  
\]
where $a \in {r_{i}(\boldsymbol{\beta},\bg)}$ indicates that route
$r_{i}(\boldsymbol{\beta},\bg)$ uses link $a$. Note also that we have assumed
existence of the partial derivatives; if not, one can replace them with
subgradients. Such partial derivatives typically do not have an exact analytical
expression and we in turn use this approximation technique; a comprehensive
discussion on this approximation can be found in \cite{Patriksson2004}. 

Similar to \cite{Zhang2016,Zhang2018}, the reasons we consider only the shortest
routes for the purpose of calculating these gradients include: $(1)$ GPS navigation
is widely-used by vehicle drivers so they tend to always select the fastest routes
between their OD pairs. $(2)$ Considering the fastest routes only significantly
simplifies the calculation of the route-choice probabilities. $(3)$ Extensive
numerical experiments show that such an approximation of the gradients performs
satisfactorily well.

\subsubsection{Directional flow derivatives with respect to parameters of the latency
function} 

To the best of our knowledge there are two main approaches \cite{TobinFriesz88,
  Patriksson2004} to calculate directional derivatives of the cost function with
respect to a perturbation $\rho$ on the cost coefficients $\bbeta$. In
\cite{TobinFriesz88} the sensitivity analysis is made with respect to the routes and
requires solving a linear system that in some cases may be difficult when dealing
with large-scale networks as pointed out in \cite{JOSEFSSON20074}. To overcome this
issue, \cite{JOSEFSSON20074} proposes a QP formulation to calculate such
derivatives. To find a solution to this QP, \cite{JOSEFSSON20074} solves a similar
problem to TAP. Therefore, although we are able to use any of these methods to
calculate $\partial x_a (\bbeta_j,\bg_j)/\partial \beta_{l}$ we prefer to use a
finite-difference approximation. This is because: $(1)$ the complexity of solving the
TAP is similar to that of the QP proposed by \cite{JOSEFSSON20074}, and $(2)$ the MSA
algorithm is an efficient algorithm that allows us to include all routes connecting
an OD pair $\bw_i$ in its route set $\mathcal{R}^{\bw_i}$. Using $\text{TAP}_a(\cdot)$
to denote the outcome of MSA for link $a$, for some small enough $\rho$ we compute
\[
\frac{\partial x_a (\bbeta_j,\bg_j)}{\partial \beta_l} \approx \frac{\text{TAP}_a(\bbeta_j+\rho
  \mb{e}_l,\bg_j) - \text{TAP}_a(\bbeta_j,\bg_j)}{\rho}, 
\]
where $\mb{e}_l$ is the $l$th unit vector.

\section{NUMERICAL EXAMPLE}
\label{sec:numericalExample}

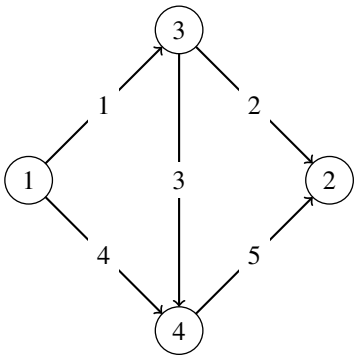
\begin{figure}[ht]
\centering
\tikzset{EdgeStyle/.append style = {->, bend left = 0} }
\thispagestyle{empty}
\begin{tikzpicture}
  \SetGraphUnit{2}
  \Vertex{1}
  \NOEA(1){3}
  \SOEA(1){4}
  \NOEA(4){2}
  \Edge[label = 1](1)(3)
  \Edge[label = 3](3)(4)
  \Edge[label = 2](3)(2)
  \tikzset{EdgeStyle/.append style = {->, bend right = 0} }
  \Edge[label = 4](1)(4)
  \Edge[label = 5](4)(2)
\end{tikzpicture}
\caption{Braess' network; we consider one OD pair from node $1$ to node
  $2$.} \label{fig:BraessNetwork} 
\end{figure}
We perform a numerical experiment to test our method. To do so, we generate
\emph{ground truth} data by choosing specific OD demands and cost functions. Then, we
solve the TAP to obtain data flows $\mathbf{x}^*$. Once we have the \emph{ground
  truth} information, we initialize our method with a feasible $f(\cdot)$ and
$\bg_0$. We aim to adjust these initial OD demands and cost functions such that the
resulting link flows $\mathbf{x}(\boldsymbol{\beta}, \bg)$ are close to the
\emph{ground truth} flows $\mathbf{x}^*$.

As an example we use the Braess network (Fig. \ref{fig:BraessNetwork}). In this
network, we generate \emph{ground truth} by considering a single OD pair which
transports $4,000$ vehicles from node $1$ to $2$. Furthermore, we consider the cost
function to be $f(x) = 1 + x$. The resulting flows when solving the TAP for this
example are: $(2080, 2080, 0, 1920, 1920)$ for links $(1,2,3,4,5)$ respectively.

Then, for solving the bilevel problem, we set an initial demand $\bg_0$ to be $5,500$
vehicles, and initial cost function equal to BPR i.e. $f(x) = 1+ 0.15x^4$, i.e.,
$\boldsymbol{\beta}_0 = (1,0,0,0,0.15,0)$. Then, we implement our model using $c=30$,
$\lambda = 10^3$, $c1=c2=5$, $\rho = 0.5$ and $n$ (polynomial degree) equal to
$5$. Notice that these parameters can be selected using cross-validation.

\begin{figure}[ht]
\centering
\includegraphics[width=0.95\columnwidth]{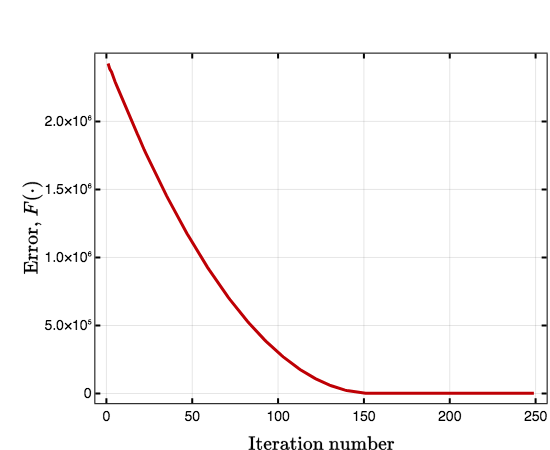}
\caption{Objective function of the Bilevel problem, i.e., $F(\boldsymbol{\beta}_j,
  \bg_j)$ as a function of the number of iterations $j$.} \label{fig:Convergance}
\end{figure}

By running experiments, we observe that the objective function of the bilevel problem
(cf. (\ref{eq:bilevel})) converges to zero (see Fig.~\ref{fig:Convergance}). However,
we also noticed that is quite sensitive to the parameters used, in particular, we
have to be careful when selecting ($c1$, $c2$) and $\lambda$ because these may cause
unboundness by violating the (IP-1) constraint set and the bilevel primal-dual gap
respectively. Moreover, note that the selection of ($c1$, $c2$) has a direct impact
on the algorithm's convergence rate.

\begin{figure}[ht]
\centering
\includegraphics[width=0.95\columnwidth]{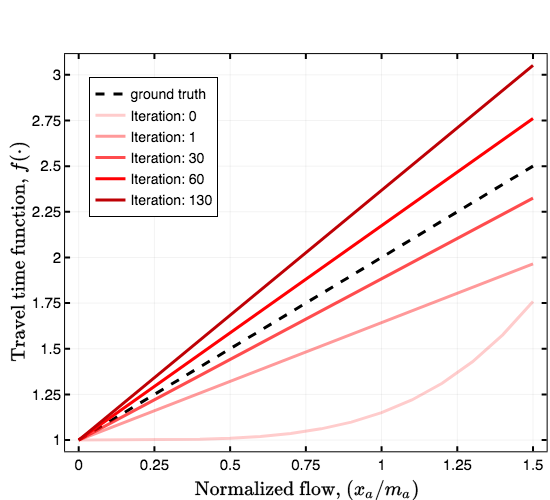}
\caption{Cost function estimators with respect to the joint iterations. In this example,  the cost function coefficient converges around iteration $j=130$.} \label{fig:CostFunctionIters}
\end{figure}

\begin{figure}[ht]
\centering
\includegraphics[width=0.95\columnwidth]{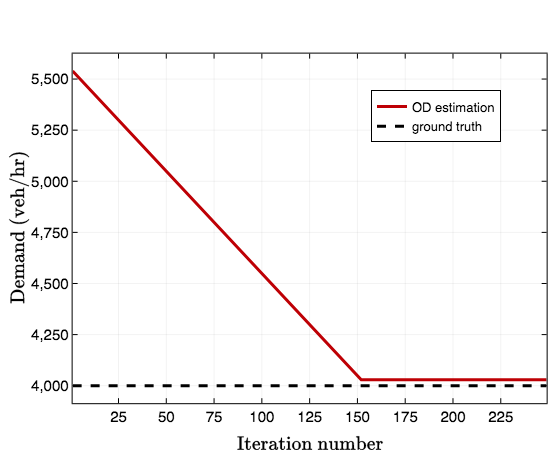}
\caption{Demand estimator for OD pair $(1,2)$ with respect to the joint iterations.}
\label{fig:demandConvergance}
\end{figure}

When solving the problem we obtain the estimated OD demand, cost function and link
flows as: $4,035$ (Fig.~\ref{fig:demandConvergance}); $f(x)=1+1.45x$
(Fig.~\ref{fig:CostFunctionIters}); and $\mathbf{x}=(2079.5, 2079.5, 0, 1950.5,
1950.5)$, respectively. This is a very good estimate of the ground truth.  Even
though the latency function is not exactly the same, it is returning similar
flows. This happens because commuters respond equally to $f(x)=1+x$ and to
$f(x)=1+1.45x$ for this particular network and conditions. We would expect the
difference between cost function estimation to decrease as we add more data samples
to the joint problem.

\section{CONCLUSION}
\label{sec:conclusion}
In this work, we were able to solve the joint problem of estimating OD demands and
cost functions in a transportation network. We approached the problem by rewriting
\eqref{eq:bilevel} with the lower-level KKT conditions
\eqref{eq:bilevel_with_kkt}. Then, we solved the problem using an iterative approach
\eqref{eq:FrankWolfe}. To be able to accomplish this, we relaxed some constraints by
allowing a small gap to exist between the primal-dual costs. Additionally, we took
care of the non-convexity of the constraints by using the previous iteration solution
and bounding these variables.

Finally, we tested our algorithm using the Braess network and concluded that our
proposed method works well in terms of reducing the objective function of the bilevel
formulation \eqref{eq:bilevel}. We performed this task by adjusting both, the OD
demand and the cost functions. It is important to keep in mind that the output of the
algorithm is sensitive to the accuracy of flow observations and to the parameters
chosen. To overcome the parameter selection issue, we suggest practitioners to use
cross-validation techniques. As future extensions of this work, we plan to implement
this algorithm in significantly larger networks and we aim at extending our framework
to multi-class transportation networks.




\bibliographystyle{IEEEtran}

\bibliography{report}

\begin{thebibliography}{10}
\providecommand{\url}[1]{#1}
\csname url@samestyle\endcsname
\providecommand{\newblock}{\relax}
\providecommand{\bibinfo}[2]{#2}
\providecommand{\BIBentrySTDinterwordspacing}{\spaceskip=0pt\relax}
\providecommand{\BIBentryALTinterwordstretchfactor}{4}
\providecommand{\BIBentryALTinterwordspacing}{\spaceskip=\fontdimen2\font plus
\BIBentryALTinterwordstretchfactor\fontdimen3\font minus
  \fontdimen4\font\relax}
\providecommand{\BIBforeignlanguage}[2]{{%
\expandafter\ifx\csname l@#1\endcsname\relax
\typeout{** WARNING: IEEEtran.bst: No hyphenation pattern has been}%
\typeout{** loaded for the language `#1'. Using the pattern for}%
\typeout{** the default language instead.}%
\else
\language=\csname l@#1\endcsname
\fi
#2}}
\providecommand{\BIBdecl}{\relax}
\BIBdecl

\bibitem{Merchant78}
D.~K. Merchant and G.~L. Nemhauser, ``A model and an algorithm for the dynamic
  traffic assignment problems,'' \emph{Transportation Science}, vol.~12, no.~3,
  pp. 183--199, 1978.

\bibitem{Patriksson1994}
M.~Patriksson, ``{The Traffic Assignment Problem: Models and Methods},''
  \emph{Annals of Physics}, vol.~54, no.~2, pp. xii, 223 p., 1994.

\bibitem{bertsimas2015data}
D.~Bertsimas, V.~Gupta, and I.~C. Paschalidis, ``Data-driven estimation in
  equilibrium using inverse optimization,'' \emph{Mathematical Programming},
  vol. 153, no.~2, pp. 595--633, 2015.

\bibitem{VanZuylen80}
H.~J. Van~Zuylen and L.~G. Willumsen, ``The most likely trip matrix estimated
  from traffic counts,'' \emph{Transportation Research Part B: Methodological},
  vol.~14, no.~3, pp. 281--293, 1980.

\bibitem{Hazelton2000}
M.~L. Hazelton, ``{Estimation of origin-destination matrices from link flows on
  uncongested networks},'' \emph{Transportation Research Part B:
  Methodological}, vol.~34, no.~7, pp. 549--566, 2000.

\bibitem{Spiess87}
H.~Spiess, ``A maximum likelihood model for estimating origin-destination
  matrices,'' \emph{Transportation Research Part B: Methodological}, vol.~21,
  no.~5, pp. 395 -- 412, 1987.

\bibitem{Daamen2014}
C.~B. Winnie~Daamen, ``Traffic simulation and data: Validation methods and
  applications,'' \emph{CRC Press}, vol. 978-1482228700, no.~1, 2014.

\bibitem{BPR}
T.~A. Manual, ``Bureau of public roads,'' \emph{US Department of Commerce},
  1964.

\bibitem{Zhang2016}
J.~Zhang, S.~Pourazarm, C.~G. Cassandras, and I.~C. Paschalidis, ``{The price
  of anarchy in transportation networks by estimating user cost functions from
  actual traffic data},'' \emph{2016 IEEE 55th Conference on Decision and
  Control, CDC 2016}, no. Cdc, pp. 789--794, 2016.

\bibitem{Zhang2018}
------, ``{The Price of Anarchy in Transportation Networks: Data-Driven
  Evaluation and Reduction Strategies},'' \emph{Proceedings of the IEEE}, vol.
  106, no.~4, 2018.

\bibitem{Yang2001}
H.~Yang, Q.~Meng, and M.~G.~H. Bell, ``Simultaneous estimation of the
  origin-destination matrices and travel-cost coefficient for congested
  networks in a stochastic user equilibrium,'' \emph{Transportation Science},
  vol.~35, no.~2, pp. 107--123, 2001.

\bibitem{Braess2005}
D.~Braess, A.~Nagurney, and T.~Wakolbinger, ``{On a Paradox of Traffic
  Planning},'' \emph{Transportation Science}, vol.~39, no.~4, pp. 446--450,
  2005.

\bibitem{Beckmann1955}
M.~J. Beckmann, C.~B. McGuire, and C.~B. Winsten, ``{Studies in the Economics
  of Transportation},'' p. 359, 1955.

\bibitem{Smith1979}
M.~J. Smith, ``{The existence, uniqueness and stability of traffic
  equilibria},'' \emph{Transportation Research Part B}, vol.~13, no.~4, pp.
  295--304, 1979.

\bibitem{Dafermos1980}
S.~Dafermos, ``{Traffic Equilibrium and Variational Inequalities},''
  \emph{Transportation Science}, vol.~14, no.~1, pp. 42--54, 1980.

\bibitem{Daganzo77}
C.~F. Daganzo and Y.~Sheffi, ``On stochastic models of traffic assignment,''
  \emph{Transportation Science}, vol.~11, no.~3, pp. 253--274, 1977.

\bibitem{Patriksson2004}
M.~Patriksson, ``Sensitivity analysis of traffic equilibria,''
  \emph{Transportation Science}, vol.~38, no.~3, pp. 258--281, Aug. 2004.

\bibitem{TobinFriesz88}
R.~Tobin and T.~Friesz, ``\BIBforeignlanguage{English (US)}{Sensitivity
  analysis for equilibrium network flow},'' \emph{\BIBforeignlanguage{English
  (US)}{Transportation Science}}, vol.~22, no.~4, pp. 242--250, 1 1988.

\bibitem{JOSEFSSON20074}
M.~Josefsson and M.~Patriksson, ``Sensitivity analysis of separable traffic
  equilibrium equilibria with application to bilevel optimization in network
  design,'' \emph{Transportation Research Part B: Methodological}, vol.~41,
  no.~1, pp. 4 -- 31, 2007.

\end{thebibliography}

\end{document}